\documentclass[12pt]{amsart}
\usepackage{amssymb}
\usepackage{amsfonts}
\usepackage{graphicx}
\usepackage{amsmath}
\usepackage{xcolor}

\newtheorem{theorem}{Theorem}

\begin{document}

\title{There is no degree independent Bombieri type inequality for non-homogeneous polynomials}
\author{ J. M. Aldaz}
\address{J.M. Aldaz: Instituto de Ciencias Matem\'aticas (CSIC-UAM-UC3M-UCM)
and Departamento de Matem\'aticas, Universidad Aut\'onoma de Madrid,
Cantoblanco 28049, Madrid, Spain.}
\email{jesus.munarriz@uam.es}

\thanks{2020 Mathematics Subject Classification: \emph{Primary: 26D05}, 
\emph{Secondary:  15A63}}
\thanks{Key words and phrases: \emph{Apolar inner product,  Bombieri's inequality}}


\maketitle

It is noted in \cite{ABR} that Bombieri's inequality can be extended to non-homogeneous polynomials just by homogeneizing them (which leads to an increase in the degrees and thus the corresponding norms) and then using the standard Bombieri inequality. Morally, Bombieri's inequality tells us that multiplying polynomials cannot lead to susbtantial cancellations, as measured by the apolar norm, when at least one of the polynomials is homogeneous. Therefore, a natural question is whether a form of the inequality, independent of the degrees, can be found when both polynomials are non-homogeneous, or in other words, whether large cancellations cannot occur in this latter case. 

\vskip .2 cm

After some 
fruitless reflection I decided to consult Gemini (flash mode) which quickly provided a one variable counterexample. In retrospect, it is a rather obvious one:
let $t > 0$ be fixed, and  let   $P_{t,m}(z)$ and $Q_{t,m}(z)$ be the $m$-th order Taylor polynomials of $e^{tz}$ and $e^{-tz}$ respectively. By choosing $m = m(t)$ sufficiently large we can make $\|P_{t,m}Q_{t,m}\|_a$ arbitrarily close to $1$
and $\|P_{t,m}\|_a\|Q_{t,m}\|_a$ arbitrarily close to $e^{t^2}$. Thus, there is no (absolute) constant $C > 0$ such that the inequality 
$$
C \|P_{t,m}Q_{t,m}\|_a \ge \|P_{t,m}\|_a\|Q_{t,m}\|_a
$$
holds for every $t > 0$ and every $m\ge 1$.

\vskip .2 cm

So this is the counterexample. Next we  fill in the details of the preceding assertions. Let us recall the definition of the apolar inner product $\langle  \cdot,\cdot\rangle_{a}$.
Denote by $\mathbb{C}[z_{1},\dots,z_{d}]$ the vector space of all
polynomials  in  $z=\left(  z_{1},\dots,z_{d}\right)$ with complex coefficients.
Let $P$ and $Q$ be polynomials of degree $N$ and $M$, respectively given by
\[
P\left(  z \right)  =\sum_{\alpha\in\mathbb{N}^{d},\left\vert
\alpha\right\vert \leq N}c_{\alpha}z ^{\alpha}\text{ and }Q\left(  z\right)
=\sum_{\alpha\in\mathbb{N}^{d},\left\vert \alpha\right\vert \leq
M}d_{\alpha}z^{\alpha},
\]
where the standard notation for multi-indices is used.

Given $P\left(  z\right)  $, we denote by
$P^{\ast}\left(  z\right)  $ the polynomial obtained from $P\left(  z\right)
$ by conjugating its coefficients, and by $P\left(  D\right)  $  the linear
differential operator obtained from $P\left(  z\right)
$ by replacing each variable $z_{j}$ with the
differential operator $\frac{\partial}{\partial z_{j}}$, for $j = 1, \dots ,d$.
The \emph{apolar inner product} $\langle  \cdot,\cdot\rangle_{a}$ on
$\mathbb{C}[z_{1},\dots,z_{d}]$ is defined by
\begin{equation}\label{apolar}
\left\langle P,Q\right\rangle_{a} :=\left( Q^*\left(
D\right)  P\right) \;(0)=\sum_{\alpha\in\mathbb{N}^{d}}\alpha!c_{\alpha
}\overline{d_{\alpha}},
\end{equation}
and the associated {\em apolar norm}, by
$
\left\Vert f\right\Vert _{a}=\sqrt{\left\langle f,f\right\rangle_{a}}.
$
Note that for monomials $z^\alpha$ and $z^\beta$, definition (\ref{apolar}) means that $\left\langle z^\alpha, z^\beta \right\rangle_{a} = 0$ if $\alpha \ne \beta$, and $\left\langle z^\alpha, z^\alpha \right\rangle_{a} = \alpha !$. Furthermore, from these relations definition (\ref{apolar}) is recovered by sesquilinearity.

\vskip .2 cm

Some computations follow:
\begin{equation*}
	P_{t,m}(z) = \sum_{k=0}^m \frac{t^k}{k!} z^k
\mbox{ \ \ \  and \ \ \ }
	Q_{t,m}(z) = \sum_{k=0}^m \frac{(-1)^k t^k}{k!} z^k.
\end{equation*}
Thus, 
$$
 \|P_{t,m}(z)\|_a\|Q_{t,m}(z)\|_a
 =
\|P_{t,m}(z)\|_a^2
=
\sum_{k=0}^m \frac{t^{2k}}{k!}.
 $$
Let us recall a well known integral representation of the apolar inner product due to  V. Bargmann, cf. \cite{Bar}.

\begin{theorem} \label{barg}
	Let $P$ and $Q$ be polynomials in $d$ complex variables. Then
	\begin{equation*}
		\left\langle P,Q\right\rangle _{a}=\frac{1}{\pi^{d}}\int_{\mathbb{R}^{d}}
		\int_{\mathbb{R}^{d}}P\left(  x+iy\right)  \overline{Q\left(  x+iy\right)
		}e^{-\left|  x\right|  ^{2}-\left|  y\right|  ^{2}}dxdy<\infty\label{Bargman},
	\end{equation*}
	where $dxdy$ is Lebesgue measure on $\mathbb{R}^{2d}$.
\end{theorem}

Since $	P_{t,m}(z)$ converges uniformly on compacta to $e^{tz}$ as $m\to \infty$,  and so does 
$Q_{t,m}(z)$ to $e^{- tz}$, it follows that the product
 $	P_{t,m}(z) \ Q_{t,m}(z)$ converges uniformly on compacta to $1$. Thus, the interchange of limits and integrals in Bargmann's representation of the apolar inner product is justified, and it follows that
  $$
 \lim_{m\to \infty} \left\langle 	P_{t,m} Q_{t,m},	P_{t,m} Q_{t,m}\right\rangle _{a}
 =
 1.
$$

As the main author of this note is Gemini (my contribution having been to check the details and make the presentation suit my taste) this ms. will not be submitted to a journal for publication. Also, I used Gemini because it is the LLM that pops out when I open the
Chrome navigator, and so far utilizing it has not cost me any money. But I have no financial interest in Gemini or any other AI. Finally, it seems that thinking first and asking an AI second, at least in this case, is not the most time efficient way to proceed.


\begin{thebibliography}{99}  

\bibitem{ABR} Aldaz J.M., Bravo, A. and  Render, H.:{\em The apolar inner product and a Bombieri-type inequality for polynomials}.  Journal of Inequalities and Applications,
Volume 2025, article number 90 (2025). 
 

	\bibitem{Bar}  Bargmann, V.:  {\em On a Hilbert space of analytic functions and an associated integral transform}.  Comm. Pure Appl. Math. \textbf{14},  (1961)  187--214. 

\end{thebibliography}
\end{document}